\documentclass{amsart}
\usepackage{amsfonts, amssymb}
\theoremstyle{plain}
 \newtheorem{thm}{Theorem}[section]
 \newtheorem{lem}[thm]{Lemma}
 \newtheorem{cor}[thm]{Corollary}

\theoremstyle{definition}
  \newtheorem{defn}{Definition}[section]

\theoremstyle{remark}
  \newtheorem{rem}{Remark}[section]

\newcommand{\ci}[2]{\cite[#1]{#2}}

\renewcommand{\c}{\curvearrowright}

\begin{document}
\title[Orbit equivalence rigidity for ergodic actions]{Orbit equivalence rigidity for ergodic actions \\ of the mapping class group}
\author{Yoshikata Kida}
\address{Mathematical Institute, Tohoku University, Sendai 980-8578, Japan}
\email{kida@math.tohoku.ac.jp}

\begin{abstract}
We establish orbit equivalence rigidity for any ergodic, essentially free and measure-preserving action on a standard Borel space with a finite positive measure of the mapping class group for a compact orientable surface with higher complexity. We prove similar rigidity results for a finite direct product of mapping class groups as well.
\end{abstract}

\maketitle

\footnote[0]{{\it Date}: June 1, 2007.}

\section{Introduction}\label{sec-int}

This paper is a continuation of \cite{kida-mer}, in which we established rigidity for the mapping class group in terms of measure equivalence. In this paper, we study ergodic standard actions of the mapping class group from the viewpoint of orbit equivalence theory, and establish several rigidity results. In this paper, by a discrete group we mean a discrete and countable group. A {\it standard} action of a discrete group means an essentially free, measure-preserving action on a standard finite measure space (i.e., a standard Borel space with a finite positive measure).

In ergodic theory, it is one of the main problems to study conjugacy classes of ergodic actions of a given group.

\begin{defn}
Consider a measure-preserving action of a discrete group $\Gamma_{i}$ on a standard finite measure space $(X_{i}, \mu_{i})$ for $i=1, 2$. The two actions are said to be {\it conjugate} if there are conull Borel subsets $X_{1}'\subset X_{1}$, $X_{2}'\subset X_{2}$, a Borel isomorphism $f\colon X_{1}'\rightarrow X_{2}'$ and an isomorphism $F \colon \Gamma_{1}\rightarrow \Gamma_{2}$ such that
\begin{enumerate}
\item[(i)] the two measures $f_{*}\mu_{1}$ and $\mu_{2}$ are equivalent;
\item[(ii)] for any $g\in \Gamma_{1}$ and any $x\in X_{1}'$, $gx$ belongs to $X_{1}'$ and the equation $f(gx)=F(g)f(x)$ holds.
\end{enumerate}
\end{defn}

On the other hand, orbit equivalence is a much weaker equivalence relation among measure-preserving actions of discrete groups on standard finite measure spaces than conjugacy.

\begin{defn}
Consider a measure-preserving action of a discrete group $\Gamma_{i}$ on a standard finite measure space $(X_{i}, \mu_{i})$ for $i=1, 2$. Then the two actions are said to be {\it weakly orbit equivalent (WOE)} if there are Borel subsets $A_{1}\subset X_{1}$ and $A_{2}\subset X_{2}$ satisfying $\Gamma_{1}A_{1}=X_{1}$ and $\Gamma_{2}A_{2}=X_{2}$ up to null sets and there is a Borel isomorphism $f\colon A_{1}\rightarrow A_{2}$ such that
\begin{enumerate}
\item[(i)] the two measures $f_{*}(\mu_{1}|_{A_{1}})$ and $\mu_{2}|_{A_{2}}$ are equivalent;
\item[(ii)] $f(\Gamma_{1}x\cap A_{1})=\Gamma_{2}f(x)\cap A_{2}$ for a.e.\ $x\in A_{1}$.
\end{enumerate}
If we can take both $A_{1}$ and $A_{2}$ to have full measure, then the two actions are said to be {\it orbit equivalent (OE)}.
\end{defn}

It is clear that two conjugate actions are OE. The study of orbit equivalence was initiated by Dye \cite{dye}, \cite{dye2}, who studied standard actions of some special amenable groups, and Ornstein and Weiss \cite{ow} concluded that any ergodic standard actions of any two infinite amenable groups are OE. More generally, Connes, Feldman, and Weiss \cite{cfw} showed that amenable discrete measured equivalence relations are hyperfinite, which implies uniqueness of ergodic amenable equivalence relations of type ${\rm II}_{1}$. On the other hand, one can easily construct a family of continuously many ergodic standard actions of $\mathbb{Z}$ which are mutually non-conjugate. These phenomena give rise to a sharp difference between orbit equivalence and conjugacy for ergodic actions of $\mathbb{Z}$. 

In contrast, based on Zimmer's pioneering work \cite{zim2}, Furman \cite{furman2} established OE superrigidity for some ergodic standard actions of a lattice in a simple Lie group of higher real rank. Given an ergodic standard action $\alpha$ of a discrete group, we say that $\alpha$ is {\it OE superrigid} if the following holds: Let $\beta$ be any ergodic standard action of a discrete group which is WOE to $\alpha$. Then $\alpha$ and $\beta$ are virtually conjugate in the following sense.

\begin{defn}
Let $\Gamma$ and $\Lambda$ be discrete groups. Suppose that they admit ergodic standard actions $\Gamma \c (X, \mu)$ and $\Lambda \c (Y, \nu)$, where $(X, \mu)$ and $(Y, \nu)$ are standard finite measure spaces. Then the two actions are said to be {\it virtually conjugate} if we can find exact sequences
\[1\rightarrow N\rightarrow \Gamma \rightarrow \Gamma_{1}\rightarrow 1,\ \ \ 
1\rightarrow M\rightarrow \Lambda \rightarrow \Lambda_{1}\rightarrow 1\]
of groups, where $N$ and $M$ are both finite, and there exist finite index subgroups $\Gamma_{2}<\Gamma_{1}$ and $\Lambda_{2}<\Lambda_{1}$ satisfying the following: Put
\[(X_{1}, \mu_{1})=(X, \mu)/N,\ \ (Y_{1}, \nu_{1})=(Y, \nu)/M\]
and consider the natural actions $\Gamma_{1}\c (X_{1}, \mu_{1})$ and $\Lambda_{1}\c (Y_{1}, \nu_{1})$.
\begin{enumerate}
\item[(i)] The action $\Gamma_{1}\c (X_{1}, \mu_{1})$ is conjugate with the induction $(X_{2}, \mu_{2})\uparrow_{\Gamma_{2}}^{\Gamma_{1}}$ from some ergodic standard action $\Gamma_{2}\c (X_{2}, \mu_{2})$. Similarly, $\Lambda_{1}\c (Y_{1}, \nu_{1})$ is conjugate with the induction $(Y_{2}, \nu_{2})\uparrow_{\Lambda_{2}}^{\Lambda_{1}}$ from some ergodic standard action $\Lambda_{2}\c (Y_{2}, \nu_{2})$ (see Definition \ref{defn-ind} for inductions).
\item[(ii)] The actions $\Gamma_{2}\c (X_{2}, \mu_{2})$ and $\Lambda_{2}\c (Y_{2}, \nu_{2})$ are conjugate.  
\end{enumerate}
\end{defn}

It is easy to see that two virtually conjugate actions are WOE. Let $\Gamma$ be a lattice in a connected simple Lie group $G$ of non-compact type with finite center and real rank at least $2$. In \cite{furman2}, Furman established OE superrigidity of some ergodic standard actions of $\Gamma$ (e.g., the standard action of $SL(n, \mathbb{Z})$ on $\mathbb{R}^{n}/\mathbb{Z}^{n}$ with $n\geq 3$). Moreover, he showed that all other ergodic standard actions of $\Gamma$ essentially come from the $\Gamma$-action on $G/\Lambda$ for some lattice $\Lambda$ in $G$, which is WOE to the $\Lambda$-action on $G/\Gamma$. Monod and Shalom \cite{ms} applied the theory of bounded cohomology to the setting of orbit equivalence, and established (slightly weaker) OE superrigidity results for irreducible standard actions of discrete groups of the form $\Gamma_{1}\times \cdots \times \Gamma_{n}$ with $n\geq 2$ and $\Gamma_{i}\in \mathcal{C}$ torsion-free for each $i$, where $\mathcal{C}$ is the class of discrete groups introduced in \cite{ms}. (A measure-preserving action of a discrete group of the form $\Gamma_{1}\times \cdots \times \Gamma_{n}$ is said to be irreducible if for every $j$, the product $\prod_{i\neq j}\Gamma_{i}$ acts ergodically.) This class $\mathcal{C}$ is huge and contains all non-elementary word-hyperbolic groups, the free product of two infinite groups (see Section 7 in \cite{ms}) and mapping class groups \cite{ham2}. Recently, Popa \cite{popa-mal2}, \cite{popa} discovered that the Bernoulli action of an infinite discrete group satisfying Kazhdan's property (T) admits stronger rigidity properties in terms of its associated von Neumann algebra. The reader is referred to \cite{vaes} for Popa's recent breakthrough rigidity results on Bernoulli actions of infinite groups satisfying property (T).

In this note, following Furman's technique in \cite{furman2} and using rigidity results due to the author developed in \cite{kida}, \cite{kida-mer}, we establish OE superrigidity for any ergodic standard action of the mapping class group. Furthermore, we show the same rigidity phenomenon for any ergodic standard action of a direct product of mapping class groups. As mentioned above, it follows from Monod and Shalom's rigidity result that irreducible standard actions of a non-trivial direct product of torsion-free finite index subgroups of mapping class groups have (slightly weaker) OE superrigidity.

Throughout the paper, we assume a surface to be connected, compact and orientable unless otherwise stated. We write $\kappa(M)=3g+p-4$ for a surface $M=M_{g, p}$ of genus $g$ and with $p$ boundary components. Let $\Gamma(M)^{\diamond}$ be the extended mapping class group of $M$, the group of isotopy classes of all diffeomorphisms of $M$. The mapping class group $\Gamma(M)$ of $M$ is the group of isotopy classes of all orientation-preserving diffeomorphisms of $M$, which is a subgroup of index $2$ in $\Gamma(M)^{\diamond}$.

\begin{thm}\label{thm-sor}
Let $n$ be a positive integer and let $M_{i}$ be a surface with  $\kappa(M_{i})>0$ for $i\in \{ 1, \ldots, n\}$. If $\Gamma$ is a finite index subgroup of $\Gamma(M_{1})^{\diamond}\times \cdots \times \Gamma(M_{n})^{\diamond}$, then any ergodic standard action of $\Gamma$ on a standard finite measure space is OE superrigid.
\end{thm}

Remark that any finite index subgroup of $\Gamma(M)^{\diamond}$ with $\kappa(M)>0$ and $M\neq M_{1, 2}, M_{2, 0}$ has no non-trivial finite normal subgroups (see 11.5 in \cite{ivanov1}). Hence, so does any finite index subgroup of $\Gamma(M_{1})^{\diamond}\times \cdots \times \Gamma(M_{n})^{\diamond}$  with $\kappa(M_{i})>0$ and $M_{i}\neq M_{1, 2}, M_{2, 0}$ for all $i$. 

Moreover, we consider the following three particular cases. Recall that a measure-preserving action of a discrete group $\Gamma$ on a measure space is said to be {\it aperiodic} if any finite index subgroup of $\Gamma$ acts ergodically on the measure space.

\begin{thm}\label{thm-sr}
Let $n$ be a positive integer and let $M_{i}$ be a surface with  $\kappa(M_{i})>0$ and $M_{i}\neq M_{1, 2}, M_{2, 0}$ for each $i\in \{ 1, \ldots, n\}$. Let $\Gamma$ be a finite index subgroup of $\Gamma(M_{1})^{\diamond}\times \cdots \times \Gamma(M_{n})^{\diamond}$ and let $\Lambda$ be a discrete group. Let $(X, \mu)$ and $(Y, \nu)$ be standard finite measure spaces. Suppose that $\Gamma$ admits an aperiodic standard action $\alpha$ on $(X, \mu)$ and that $\Lambda$ admits an ergodic standard action $\beta$ on $(Y, \nu)$. Then the following two assertions hold:
\begin{enumerate}
\item[(i)] Suppose that $\Lambda$ has no non-trivial finite normal subgroups and that the action $\beta$ is also aperiodic. If the two actions $\alpha$ and $\beta$ are WOE, then they are conjugate.
\item[(ii)] If the two actions $\alpha$ and $\beta$ are OE, then they are conjugate.
\end{enumerate}
\end{thm}

\begin{thm}\label{thm-ssr}
Let $\Gamma$ be as in Theorem \ref{thm-sr}. Let $(X, \mu)$ and $(Y, \nu)$ be standard finite measure spaces. Suppose that $\Gamma$ admits an aperiodic standard action on $(X, \mu)$ and an ergodic standard action on $(Y, \nu)$. If the two actions are WOE, then they are conjugate.
\end{thm}

\begin{thm}\label{thm-ssor}
Let $M$ be a surface with $\kappa(M)>0$ and $M\neq M_{1, 2}, M_{2, 0}$. Let $\Gamma$ be a normal subgroup of finite index in $\Gamma(M)^{\diamond}$. If two ergodic standard actions of $\Gamma$ are WOE, then they are conjugate.
\end{thm}

The same rigidity property as in Theorem \ref{thm-ssor} is satisfied for certain subgroups of a direct product of  mapping class groups as well (see Corollary \ref{cor-woe-iso}).

\begin{rem}
Thanks to these rigidity results, we can construct a new example of an ergodic equivalence relation of type ${\rm II}_{1}$ which cannot arise from any standard action of a discrete group (see Corollary \ref{cor-eq-rel}). The first construction of such an equivalence relation is due to Furman \cite{furman2}, and it solved a longstanding problem formulated by Feldman and Moore \cite{fm}. 
\end{rem}

In the subsequent papers \cite{kida-ber}, \cite{kida-out}, we give applications of the above rigidity results. In \cite{kida-ber}, several classification results of generalized Bernoulli actions of mapping class groups up to OE are given. In \cite{kida-out}, we study the outer automorphism groups of discrete measured equivalence relations arising from ergodic standard actions of mapping class groups. We explicitly compute the outer automorphism groups for some special actions.

\vspace{1em}

\noindent {\bf Acknowledgements.} The author is grateful to Professor Ursula Hamenst\"adt for reading the first draft of this paper very carefully and giving many valuable suggestions. This work was done during the stay at the Max Planck Institute for Mathematics in Bonn, and the author thanks them for warm hospitality.


\section{Orbit equivalence rigidity}

In this section, we prove orbit equivalence rigidity results. For that, we fix the following notation and recall some known results. For a surface $M=M_{g, p}$ of genus $g$ and with $p$ boundary components, put $\kappa(M)=3g+p-4$. Let $M$ be a surface with $\kappa(M)>0$. Let $C=C(M)$ be the curve complex, which is a simplicial complex defined as follows. Let $V(C)=V(C(M))$ be the set of vertices of the curve complex of $M$, that is, the set of all non-trivial isotopy classes of non-peripheral simple closed curves on $M$. Let $S(M)$ be the set of its simplices, that is, the set of all non-empty finite subsets of $V(C)$ which can be realized disjointly on $M$ at the same time. Let ${\rm Aut}(C)$ be the automorphism group of the simplicial complex $C$. Then there is a natural homomorphism $\pi \colon \Gamma(M)^{\diamond}\rightarrow {\rm Aut}(C)$. It is a natural question whether $\pi$ is an isomorphism or not. The following theorem says that $\pi$ is in fact an isomorphism for almost all surfaces $M$. In \cite{ivanov3}, Ivanov sketched a proof of this statement for surfaces of genus at least 2, and Korkmaz \cite{korkmaz} gave a proof for some surfaces of genus less than 2. Luo \cite{luo} suggested another approach for this question, which does not distinguish the cases of surfaces of higher and lower genus, and finally concluded the following  

\begin{thm}\label{thm-cc-auto}
Let $M$ be a surface with $\kappa(M)>0$.
\begin{enumerate}
\item[(i)] If $M$ is neither $M_{1, 2}$ nor $M_{2, 0}$, then $\pi$ is an isomorphism.
\item[(ii)] If $M=M_{1, 2}$, then the image of $\pi$ is a subgroup of ${\rm Aut}(C)$ with index $5$ and $\ker(\pi)$ is the subgroup generated by a hyperelliptic involution, which is isomorphic to $\mathbb{Z}/2\mathbb{Z}$.  
\item[(iii)] If $M=M_{2, 0}$, then $\pi$ is surjective and $\ker(\pi)$ is the subgroup generated by a hyperelliptic involution, which is isomorphic to $\mathbb{Z}/2\mathbb{Z}$.
\item[(iv)] The two curve complexes $C(M_{0, 5})$, $C(M_{1, 2})$ (resp. $C(M_{0, 6})$, $C(M_{2, 0})$) are isomorphic as a simplicial complex.
\end{enumerate}
\end{thm}

Let $n$ be a positive integer and let $M_{i}$ be a surface with $\kappa(M_{i})>0$ and $M_{i}\neq M_{1, 2}, M_{2, 0}$ for $i\in \{ 1,\ldots, n\}$. Put $G=\Gamma(M_{1})^{\diamond}\times \cdots \times \Gamma(M_{n})^{\diamond}$. Suppose that there are a bijection $t$ on the set $\{ 1,\ldots, n\}$ and an isotopy class $\varphi_{i}$ of a diffeomorphism $M_{t(i)}\rightarrow M_{i}$. Then one can define an automorphism $\pi_{\varphi}\colon G\rightarrow G$ by
\[\pi_{\varphi}(\gamma)=(\varphi_{1}\gamma_{t(1)}\varphi_{1}^{-1}, \ldots, \varphi_{n}\gamma_{t(n)}\varphi_{n}^{-1})\]
for $\gamma =(\gamma_{1}, \ldots, \gamma_{n})\in G$. The following theorem says that all injective homomorphisms from a finite index subgroup of $G$ onto a finite index subgroup of $G$ are of the above form.  

\begin{thm}[\ci{Corollary 7.3}{kida-mer}]\label{thm-inj}
Let $\Gamma$ be a finite index subgroup of $G$. Suppose that we have an injective homomorphism $\tau \colon \Gamma \rightarrow G$ with the index $[G:\tau(\Gamma)]$ finite. Then we can find a bijection $t$ on the set $\{ 1,\ldots, n\}$ and an isotopy class $\varphi_{i}$ of a diffeomorphism $M_{t(i)}\rightarrow M_{i}$ for each $i$ such that
\[\tau(\gamma)=(\varphi_{1}\gamma_{t(1)}\varphi_{1}^{-1}, \ldots, \varphi_{n}\gamma_{t(n)}\varphi_{n}^{-1})\]
for any $\gamma =(\gamma_{1}, \ldots, \gamma_{n})\in \Gamma$. In particular, $\tau$ can be extended to an automorphism of $G$.
\end{thm}

Next, we give some elementary terminologies. If $H$, $\Lambda_{1}$, $\Lambda_{2}$ are discrete groups and $\tau_{i}\colon \Lambda_{i}\rightarrow H$ is a homomorphism for $i=1, 2$, then we denote by $(H, \tau_{1}, \tau_{2})$ the Borel space $H$ equipped with the $(\Lambda_{1}\times \Lambda_{2})$-action defined by
\[(\lambda_{1}, \lambda_{2})h=\tau_{1}(\lambda_{1})h\tau_{2}(\lambda_{2})^{-1}\]
for $h\in H$ and $(\lambda_{1}, \lambda_{2})\in \Lambda_{1}\times \Lambda_{2}$.

We shall recall a ME coupling of discrete groups. Let $\Gamma$ and $\Lambda$ be discrete groups. Suppose that $\Gamma$ and $\Lambda$ are measure equivalent, that is, there exists a measure-preserving action of $\Gamma \times \Lambda$ on a standard Borel space $(\Omega, m)$ with a $\sigma$-finite positive measure such that both of the actions $\Gamma (\simeq \Gamma \times \{ e\})\c \Omega$ and $\Lambda (\simeq \{ e\} \times \Lambda)\c \Omega$ are essentially free and have a fundamental domain of finite measure. The space $(\Omega, m)$ (equipped with the $(\Gamma \times \Lambda)$-action) is then called a {\it ME coupling} of $\Gamma$ and $\Lambda$. The reader is referred to Section 2 in \cite{furman1} and Section 3 in \cite{furman2} for fundamental properties of ME couplings.

As a fundamental fact, a ME coupling $(\Omega, m)$ of $\Gamma$ and $\Lambda$ gives rise to WOE between the actions $\Gamma \c \Omega /\Lambda$ and $\Lambda \c \Omega /\Gamma$. One can naturally identify a fundamental domain $X\subset \Omega$ for the $\Lambda$-action on $\Omega$ with $\Omega /\Lambda$ as a measure space, which induces a $\Gamma$-action on $X$. We can say the same thing for a fundamental domain $Y\subset \Omega$ for the $\Gamma$-action on $\Omega$. Conversely, the following theorem says that one can construct the corresponding ME coupling of $\Gamma$ and $\Lambda$ from WOE between ergodic standard actions of $\Gamma$ and $\Lambda$.

\begin{thm}[\ci{Theorem 3.3}{furman2}]\label{thm-furman}
Let $\Gamma$, $\Lambda$ be discrete groups and let $\Gamma \c (X, \mu)$, $\Lambda \c (Y, \nu)$ be ergodic standard actions with $\mu(X)=\nu(Y)$. Suppose that there are Borel subsets $A\subset X$, $B\subset Y$ with positive measure and a Borel isomorphism $f\colon A\rightarrow B$ such that
\begin{enumerate}
\item[(i)] the two measures $f_{*}(\mu|_{A})$ and $\nu|_{B}$ are equivalent;
\item[(ii)] $f(\Gamma x\cap A)=\Lambda f(x)\cap B$ for a.e.\ $x\in A$.
\end{enumerate}
Then there exists a ME coupling $(\Omega, m)$ of $\Gamma$ and $\Lambda$ such that
\begin{enumerate}
\item[(a)] the induced action $\Gamma \c \Omega /\Lambda$ is conjugate with $\Gamma \c X$, where the associated automorphism of $\Gamma$ is given by the identity. We can say the same thing for the actions $\Lambda \c \Omega /\Gamma$ and $\Lambda \c (Y, \nu)$;
\item[(b)] $X$ and $Y$ can be realized as fundamental domains $\bar{X}, \bar{Y}\subset \Omega$, respectively, which satisfy the equation $\mu(A)m(\bar{Y})=\nu(B)m(\bar{X})$.
\end{enumerate}
\end{thm}

In what follows, we prove the main theorems stated in Section \ref{sec-int}.

\vspace{1em}

\noindent {\it Proof of Theorem \ref{thm-sor}.} This proof heavily relies on Fumran's proof of Theorem A in \cite{furman2}. Let $n$ be a positive integer and let $M_{i}$ be a surface with $\kappa(M_{i})>0$ for $i\in \{ 1, \ldots, n\}$. Let $\Gamma$ be a finite index subgroup of
\[G=\Gamma(M_{1})^{\diamond}\times \cdots \times \Gamma(M_{n})^{\diamond}.\]
We may assume that $M_{i}\neq M_{1, 2}, M_{2, 0}$ for each $i$ by Theorem \ref{thm-cc-auto}. We identify $\Gamma(M_{i})^{\diamond}$ and ${\rm Aut}(C(M_{i}))$ via the natural isomorphism $\pi$. Let $\Lambda$ be a discrete group. Suppose that ergodic standard actions $\Gamma \c (X, \mu)$ and $\Lambda  \c (Y, \nu)$ are WOE. It follows from Theorem \ref{thm-furman} that we can construct a ME coupling $(\Omega, m)$ of $\Gamma$ and $\Lambda$ such that $\Gamma \c X$ and $\Gamma \c \Omega /\Lambda$ (resp. $\Lambda \c Y$ and $\Lambda \c \Omega /\Gamma$) are conjugate, where the associated automorphism of $\Gamma$ (resp. $\Lambda$) can be taken to be the identity. By Theorem 1.3 in \cite{kida-mer}, there exists a homomorphism $\rho \colon \Lambda \rightarrow G$ such that $\ker \rho$ and the index $[G: \rho(\Lambda)]$ are both finite. Put 
\[\Gamma_{1}=\Gamma,\ \ \Lambda_{1}=\rho(\Lambda),\ \ \Omega_{1}=\Omega /(\{ e\} \times {\rm ker}\rho ).\]
Let $m_{1}$ be the image of the measure $m$ via the quotient map $\Omega \rightarrow \Omega_{1}$. Then we have the natural $(\Gamma_{1}\times \Lambda_{1})$-action on $(\Omega_{1}, m_{1})$, and it is a ME coupling of $\Gamma_{1}$ and $\Lambda_{1}$. It follows from Corollary 7.2 in \cite{kida-mer} that we can find
\begin{enumerate}
\item[(a)] a bijection $t$ on the set $\{1,\ldots, n\}$;
\item[(b)] an isotopy class $\varphi_{i}$ of a diffeomorphism $M_{t(i)}\rightarrow M_{i}$ for each $i\in \{ 1,\ldots, n\}$;  
\item[(c)] an essentially unique almost $(\Gamma_{1}\times \Lambda_{1})$-equivariant Borel map $\Phi \colon \Omega_{1}\rightarrow (G, {\rm id}, \pi_{\varphi})$, that is,
\[\Phi((\gamma, \lambda)\omega)=\gamma \Phi(\omega)\pi_{\varphi}(\lambda)^{-1}\]
for any $\gamma \in \Gamma_{1}$, $\lambda \in \Lambda_{1}$ and a.e.\ $\omega \in \Omega_{1}$.
\end{enumerate}
Here, $\pi_{\varphi}\colon G \rightarrow G$ is the automorphism defined by
\[\pi_{\varphi}(\gamma)=(\varphi_{1}\gamma_{t(1)}\varphi_{1}^{-1}, \ldots, \varphi_{n}\gamma_{t(n)}\varphi_{n}^{-1})\]
for $\gamma =(\gamma_{1}, \ldots, \gamma_{n})\in G$. 

Put $m_{0}=\Phi_{*}m_{1}$ and let $Y_{0}, X_{0}\subset G$ be fundamental domains of $\Gamma_{1}$-, $\Lambda_{1}$-actions on the measure space $(G, {\rm id}, \pi_{\varphi})$ with the measure $m_{0}$, respectively. We can choose $Y_{0}, X_{0}$ so that $m_{0}(X_{0}\cap Y_{0})>0$. Moreover, we may assume that $m_{0}(\{ g\})>0$ for any $g\in X_{0}\cup Y_{0}$. Let $g_{0}\in X_{0}\cap Y_{0}$ and put 
\[\Gamma_{2}=\Gamma_{1}\cap g_{0}\pi_{\varphi}(\Lambda_{1})g_{0}^{-1},\ \ \Lambda_{2}=\pi_{\varphi}^{-1}(g_{0}^{-1}\Gamma_{2}g_{0}).\]
Then $\Gamma_{2}$ and $\Lambda_{2}$ are subgroups of finite index in $\Gamma_{1}$ and $\Lambda_{1}$, respectively. Note that if $(G, {\rm id}, \pi_{\varphi})$ equipped with the measure $m_{0}$ is viewed as a ME coupling of $\Gamma_{1}$ and $\Lambda_{1}$, then $\Gamma_{2}$ is the stabilizer of $g_{0}$ for the induced action $\Gamma_{1}\c X_{0}$. Similarly, $\Lambda_{2}$ is the stabilizer of $g_{0}$ for the induced action $\Lambda_{1}\c Y_{0}$. Put 
\[X_{1}=\Phi^{-1}(X_{0}),\ \ Y_{1}=\Phi^{-1}(Y_{0}),\]
which are fundamental domains of the $\Lambda_{1}$-, $\Gamma_{1}$-actions on $(\Omega_{1}, m_{1})$, respectively. Since both of the actions $\Gamma \c X$ and $\Lambda \c Y$ are ergodic, so are both of the actions $\Gamma_{1}\c X_{1}$ and $\Lambda_{1}\c Y_{1}$. Therefore, both of the actions $\Gamma_{1}\c X_{0}$ and $\Lambda_{1}\c Y_{0}$ are transitive because $\Phi$ induces a $\Gamma_{1}$-equivariant Borel map $X_{1}\rightarrow X_{0}$ and a $\Lambda_{1}$-equivariant Borel map $Y_{1}\rightarrow Y_{0}$. Put 
\[X_{2}=Y_{2}=\Phi^{-1}(\{ g_{0}\})\subset X_{1}\cap Y_{1}.\]

\begin{lem}\label{lem-inv-ind}
In the above notation,
\begin{enumerate}
\item[(i)] if we consider the induced action $\Gamma_{1}\c X_{1}$, then $X_{2}$ is $\Gamma_{2}$-invariant. Similarly, $Y_{2}$ is $\Lambda_{2}$-invariant for the induced action $\Lambda_{1}\c Y_{1}$. 
\item[(ii)] the action $\Gamma_{1}\c X_{1}$ is conjugate with the induction $X_{2}\uparrow_{\Gamma_{2}}^{\Gamma_{1}}$ from $\Gamma_{2}\c X_{2}$. Similarly, $\Lambda_{1}\c Y_{1}$ is conjugate with the induction $Y_{2}\uparrow_{\Lambda_{2}}^{\Lambda_{1}}$ from $\Lambda_{2}\c Y_{2}$. 
\end{enumerate}
\end{lem}

Assertion (i) is clear. We shall recall an induction of a group action.

\begin{defn}\label{defn-ind}
Let $\Delta$ be a discrete group and let $\Delta_{0}$ be its subgroup. Suppose that $\Delta_{0}$ admits a measure-preserving action on a measure space $(Z, \zeta)$. Define a $\Delta \times \Delta_{0}$-action on $Z\times \Delta$ by
\[(\delta, \delta_{0})(z, \delta')=(\delta_{0}z, \delta \delta'\delta_{0}^{-1})\]
for $\delta, \delta'\in \Delta$, $\delta_{0}\in \Delta_{0}$ and $z\in Z$. Then the {\it induction} from the action $\Delta_{0}\c Z$ is defined to be the natural action of $\Delta$ on the quotient space of $Z\times \Delta$ by the action of $\{ e\} \times \Delta_{0}$. We denote this action by $(Z, \zeta)\uparrow_{\Delta_{0}}^{\Delta}$.
\end{defn}

\begin{lem}\label{lem-delta}
Let $\Delta$ be a discrete group and suppose that $\Delta$ admits a measure-preserving action on a measure space $(Z, \zeta)$. Let $K$ be a countable set on which $\Delta$ acts transitively, and let $\Delta_{0}$ be the stabilizer of some $k_{0}\in K$. Suppose that there is a $\Delta$-equivariant Borel map $\Psi \colon Z\rightarrow K$. Then the action $\Delta \c Z$ is conjugate with the induction $Z_{0}\uparrow_{\Delta_{0}}^{\Delta}$ from $\Delta_{0}\c Z_{0}$, where $Z_{0}=\Psi^{-1}(k_{0})$.    
\end{lem}

\begin{proof}
The Borel map
\[f\colon Z_{0}\times \Delta \rightarrow Z,\ \ (z, \delta)\mapsto \delta z\]
satisfies that
\[f(z, \delta \delta')=\delta \delta'z=\delta f(z, \delta')\]
for $z\in Z_{0}$ and $\delta, \delta'\in \Delta$, and
\[f(\delta_{0}z, \delta'\delta_{0}^{-1})=\delta'z=f(z, \delta')\] 
for $z\in Z_{0}$, $\delta'\in \Delta$ and $\delta_{0}\in \Delta_{0}$. Therefore, $f$ induces a $\Delta$-equivariant Borel map from the quotient space of $Z_{0}\times \Delta$ by the $\Delta_{0}$-action into $Z$, which is a Borel isomorphism. 
\end{proof}

Return to the proof of Theorem \ref{thm-sor}. Assertion (ii) in Lemma \ref{lem-inv-ind} follows from Lemma \ref{lem-delta}. For the proof of Theorem \ref{thm-sor}, it is enough to show that the actions $\Gamma_{2}\c X_{2}$ and $\Lambda_{2}\c Y_{2}$ are conjugate. Note that $\Gamma_{2}g_{0}=g_{0}\pi_{\varphi}(\Lambda_{2})$ by the definition of $\Gamma_{2}$ and $\Lambda_{2}$. It follows that $\Gamma_{2}Y_{2}=\Lambda_{2}X_{2}\subset \Omega_{1}$. Let $\Omega_{2}$ denote this subspace, which is a ME coupling of $\Gamma_{2}$ and $\Lambda_{2}$, and $Y_{2}$, $X_{2}$ are fundamental domains for the $\Gamma_{2}$-, $\Lambda_{2}$-actions on $\Omega_{2}$, respectively. To distinguish from the actions on $\Omega_{2}$, we denote the induced actions $\Gamma_{2}\c X_{2}$, $\Lambda_{2}\c Y_{2}$ by $\gamma \cdot x$, $\lambda \cdot y$, using a dot. We can define a Borel map
\[\alpha \colon \Gamma_{2}\times X_{2}\rightarrow \Lambda_{2} {\rm \ \ so \ that\ \ } \gamma \cdot x=\alpha(\gamma, x)\gamma x\]
for $x\in X_{2}$ and $\gamma \in \Gamma_{2}$ because $X_{2}$ is a fundamental domain for the action $\Lambda_{2}\c \Omega_{2}$. Let $x\in X_{2}$ and $\gamma \in \Gamma_{2}$. Since $\gamma \cdot x\in X_{2}$, we see that
\[g_{0}=\Phi(\gamma \cdot x)=\Phi(\alpha(\gamma, x)\gamma x)=\gamma \Phi(x)\pi_{\varphi}(\alpha(\gamma, x))^{-1}=\gamma g_{0}\pi_{\varphi}(\alpha(\gamma, x))^{-1},\]   
and thus $\alpha(\gamma, x)=\pi_{\varphi}^{-1}(g_{0}^{-1}\gamma g_{0})$. On the other hand, if we view $x$ to be a point of $Y_{2}$, then
\[\pi_{\varphi}^{-1}(g_{0}^{-1}\gamma g_{0})\cdot x =\gamma \pi_{\varphi}^{-1}(g_{0}^{-1}\gamma g_{0})x=\gamma \cdot x,\]
where the first equality holds because the second term is in $Y_{2}$. Therefore, the isomorphism 
\[F\colon \Gamma_{2}\rightarrow \Lambda_{2},\ \ \gamma \mapsto \pi_{\varphi}^{-1}(g_{0}^{-1}\gamma g_{0})\]
and the identity $X_{2}\rightarrow Y_{2}$ give a conjugacy between the actions $\Gamma_{2}\c X_{2}$ and $\Lambda_{2}\c Y_{2}$. \hfill $\square$

\vspace{1em}

The following corollary is an immediate consequence of Theorem \ref{thm-sor}. For Assertion (ii), it is enough to note that $\Gamma$ has no non-trivial finite normal subgroups (see 11.5 in \cite{ivanov1}).

\begin{cor}\label{cor-eq-rel}
Let $\Gamma$ be as in Theorem \ref{thm-sor} and suppose that we have an ergodic standard action of $\Gamma$ on a standard finite measure space $(X, \mu)$ with $\mu(X)=1$. We denote by 
\[\mathcal{R}=\{ (\gamma x, x)\in X\times X: \gamma \in \Gamma, x\in X\}\]
the induced discrete measured equivalence relation. Then the following assertions hold:
\begin{enumerate}
\item[(i)] Let $A\subset X$ be a Borel subset with $\mu(A)$ irrational. Then the restricted relation $\mathcal{R}\cap (A\times A)$ cannot arise from any standard action of a discrete group.
\item[(ii)] Suppose that $M_{i}\neq M_{1, 2}, M_{2, 0}$ for all $i$ and that the action $\Gamma \c (X, \mu)$ is aperiodic. Then the same conclusion as in Assertion {\rm (i)} holds for any Borel subset $A\subset X$ with $0<\mu(A)<1$.
\end{enumerate} 
\end{cor}

Next, we consider special cases of Theorem \ref{thm-sor}. Let $n$ be a positive integer and let $M_{i}$ be a surface with  $\kappa(M_{i})>0$ and $M_{i}\neq M_{1, 2}, M_{2, 0}$ for each $i\in \{ 1, \ldots, n\}$. Let $\Gamma$ be a finite index subgroup of $G=\Gamma(M_{1})^{\diamond}\times \cdots \times \Gamma(M_{n})^{\diamond}$. Let $\Lambda$ be a discrete group. Let $(X, \mu)$ and $(Y, \nu)$ be standard finite measure spaces. Suppose that $\Gamma$ admits an aperiodic standard action $\alpha$ on $(X, \mu)$ and that $\Lambda$ admits an ergodic standard action $\beta$ on $(Y, \nu)$. Theorem \ref{thm-sr} (i) is then an immediate consequence of Theorem \ref{thm-sor}.

\vspace{1em}
\noindent {\it Proof of Theorem \ref{thm-sr} (ii).} We use the notation in the proof of Theorem \ref{thm-sor}. We may assume that $\mu(X)=\nu(Y)$. Since the actions $\Gamma \c X$ and $\Lambda \c Y$ are OE, we can construct a ME coupling $(\Omega, m)$ such that $Y, X\subset \Omega$ are identified with fundamental domains of the $\Gamma$-, $\Lambda$-actions on $\Omega$, respectively, and $m(X)=m(Y)$ (see Theorem \ref{thm-furman}). As in the proof of Theorem \ref{thm-sor}, put
\[\Omega_{1}=\Omega /(\{ e\} \times {\rm ker}\rho)\]
and let $m_{1}$ be the image of the measure $m$ via the quotient map $\Omega \rightarrow \Omega_{1}$. Then $m_{1}(X_{1})=m(X)\times |{\rm ker}\rho|$ and $m_{1}(Y_{1})=m(Y)$, where $|{\rm ker}\rho|$ denotes the cardinality of ${\rm ker}\rho$. Since the action $\Gamma \c X$ is aperiodic, we see that $X_{1}=X_{2}$ and $\Gamma_{1}=\Gamma_{2}$. Moreover, $m_{1}(Y_{1})=[\Lambda_{1}:\Lambda_{2}]m_{1}(Y_{2})$. Since $X_{2}=Y_{2}$, we see that ${\rm ker}\rho$ is trivial and $\Lambda_{1}=\Lambda_{2}$, which implies that the two actions $\Lambda \c Y$ and $\Lambda_{2}\c Y_{2}$ are conjugate. \hfill $\square$

\vspace{1em}
\noindent {\it Proof of Theorem \ref{thm-ssr}.} In the proof of Theorem \ref{thm-sor}, put $\Lambda =\Gamma$. Then the homomorphism $\rho \colon \Lambda \rightarrow \Lambda_{1}$ is an isomorphism because $\Lambda$ has no non-trivial finite normal subgroups (see 11.5 in \cite{ivanov1}). Therefore, we may assume that
\[\Omega =\Omega_{1},\ \ X=X_{1},\ \  Y=Y_{1}.\]
By Theorem \ref{thm-inj}, $\rho$ is the restriction of an automorphism of $G$. Thus, $[G:\Lambda]=[G:\Lambda_{1}]$. Since the action $\Gamma \c X$ is aperiodic, we see that $X_{1}=X_{2}$ and $\Gamma_{1}=\Gamma_{2}$. Note that the isomorphism $\Gamma_{2}\ni \gamma \mapsto \pi_{\varphi}^{-1}(g_{0}^{-1}\gamma g_{0})\in \Lambda_{2}$ is the restriction of an automorphism of $G$. Thus, $[G:\Gamma_{2}]=[G:\Lambda_{2}]$ and
\[[\Lambda_{1}:\Lambda_{2}]=\frac{[G:\Lambda_{2}]}{[G:\Lambda_{1}]}=\frac{[G:\Gamma_{2}]}{[G:\Lambda]}=\frac{[G:\Gamma]}{[G:\Gamma]}=1,\]
which implies that $\Lambda_{1}=\Lambda_{2}$ and $Y_{1}=Y_{2}$. \hfill $\square$

\begin{cor}\label{cor-woe-iso}
Let $n$ be a positive integer and let $M_{i}$ be a surface with $\kappa(M_{i})>0$ and $M_{i}\neq M_{1, 2}, M_{2, 0}$ for $i\in \{ 1,\ldots, n\}$. Let $\Gamma$ be a normal subgroup of finite index in $\Gamma(M_{1})^{\diamond}\times \cdots \times \Gamma(M_{n})^{\diamond}$ such that if $t$ is a bijection on the set $\{ 1, \ldots, n\}$ and $\varphi_{i}$ is an isotopy class of a diffeomorphism $M_{t(i)}\rightarrow M_{i}$ for $i\in \{ 1, \ldots, n\}$, then $\pi_{\varphi}(\Gamma)=\Gamma$ (see the sentence before Theorem \ref{thm-inj} for the definition of $\pi_{\varphi}$). Suppose that we have two ergodic standard actions $\Gamma \c (X, \mu)$ and $\Gamma \c (Y, \nu)$. If the two actions are WOE, then they are conjugate.
\end{cor}

In the notation of Corollary \ref{cor-woe-iso}, if $\Gamma$ is of the form $\Gamma_{1}\times \cdots \times \Gamma_{n}$, where $\Gamma_{i}$ is a normal subgroup of finite index in $\Gamma(M_{i})^{\diamond}$ and if $\Gamma_{i}$ and $\Gamma_{j}$ are isomorphic for any $i$, $j$ with $M_{i}$, $M_{j}$ diffeomorphic, then $\Gamma$ satisfies the hypothesis of Corollary \ref{cor-woe-iso} (see Theorem \ref{thm-inj}). In particular, if $n=1$ and $\Gamma$ is a normal subgroup of finite index in $\Gamma(M)^{\diamond}$, then the hypothesis is satisfied. Theorem \ref{thm-ssor} follows from Corollary \ref{cor-woe-iso}.

\vspace{1em}

\noindent {\it Proof of Corollary \ref{cor-woe-iso}.}
In the proof of Theorem \ref{thm-sor}, put $\Lambda =\Gamma$. Recall that $\Gamma$ has no non-trivial finite normal subgroups (see 11.5 in \cite{ivanov1}). Thus, the homomorphism $\rho \colon \Lambda \rightarrow G$ is injective. It follows from Theorem \ref{thm-inj} that $\rho =\pi_{\psi}$ for some bijection $s$ on $\{ 1, \ldots, n\}$ and an isotopy class $\psi_{i}$ of a diffeomorphism $M_{s(i)}\rightarrow M_{i}$. In the same notation as in the proof of Theorem \ref{thm-sor}, the assumption of $\Gamma$ implies that 
\[\Gamma=\Gamma_{1}=\Gamma_{2}=\Lambda_{2}=\Lambda_{1}=\Lambda\] 
by the definition of $\Gamma_{2}$ and $\Lambda_{2}$. Therefore, $X=X_{1}=X_{2}$ and $Y=Y_{1}=Y_{2}$. \hfill $\Box$

\begin{cor}\label{cor-fund}
Let $\Gamma$ be a discrete group as in Theorem \ref{thm-sor}. Suppose that $\Gamma$ admits an ergodic standard action on a standard finite measure space $(X, \mu)$. Then the fundamental group of the discrete measured equivalence relation arising from the action is trivial.
\end{cor}

\begin{proof}
Let
\[\mathcal{R}=\{ (\gamma x, x)\in X\times X: \gamma \in \Gamma, x\in X\}\]
be the discrete measured equivalence relation associated with the ergodic standard action $\Gamma \c (X, \mu)$. Recall that the fundamental group of $\mathcal{R}$ is defined to be the subgroup of the multiplicative group $\mathbb{R}_{+}^{*}$ of positive real numbers generated by $t\in \mathbb{R}_{+}^{*}$ such that for some/any Borel subset $A$ of $X$ with $\mu(A)/\mu(X)=t$, there exists an isomorphism between $\mathcal{R}$ and the restricted relation $\mathcal{R}\cap (A\times A)$. Note that for a discrete measured equivalence relation $\mathcal{S}$ of type ${\rm II}_{1}$ on a standard finite measure space $(Z, \zeta)$, and a Borel subset $B\subset Z$ of positive measure, the fundamental groups of $\mathcal{S}$ and the restricted relation $\mathcal{S}\cap (B\times B)$ are isomorphic. Therefore, for the proof of the corollary, we may assume that $M_{i}\neq M_{1, 2}, M_{2, 0}$ for all $i\in \{ 1,\ldots, n\}$ by Theorem \ref{thm-cc-auto}.

Let $A\subset X$ be a Borel subset such that $\mathcal{R}$ and $\mathcal{R}\cap (A\times A)$ are isomorphic. Put $\Lambda =\Gamma$ and let $\Lambda \c (Y, \nu)$ be a copy of $\Gamma \c (X, \mu)$. Then we can construct a ME coupling $(\Omega, m)$ of $\Gamma$ and $\Lambda$ such that $Y, X\subset \Omega$ can be identified with fundamental domains of the $\Gamma$-, $\Lambda$-actions on $\Omega$, respectively, and $\mu(A)m(Y)=\nu(Y)m(X)=\mu(X)m(X)$ (see Theorem \ref{thm-furman}). In what follows, we use the notation in the proof of Theorem \ref{thm-sor}. Since $\Lambda$ has no non-trivial finite normal subgroups, $\rho \colon \Lambda \rightarrow \Lambda_{1}$ is an isomorphism, and $m(X)=m_{1}(X_{1})$, $m(Y)=m_{1}(Y_{1})$. Since the isomorphism $\rho$ is the restriction of an automorphism of $G$ by Theorem \ref{thm-inj}, we see that
\[[G: \Gamma_{1}]=[G:\Gamma]=[G: \Lambda]=[G:\Lambda_{1}].\]
Note that
\[m_{1}(X_{1})/m_{1}(X_{2})=[\Gamma_{1}: \Gamma_{2}],\ \ m_{1}(Y_{1})/m_{1}(Y_{2})=[\Lambda_{1}:\Lambda_{2}].\]  
Since $\Gamma_{2}$ and $\Lambda_{2}$ are isomorphic via the isomorphism $\Gamma_{2}\ni \gamma \mapsto \pi_{\varphi}^{-1}(g_{0}^{-1}\gamma g_{0})\in \Lambda_{2}$, which is the restriction of an automorphism of $G$, we obtain the equation
\[[G:\Gamma_{2}]=[G:\Lambda_{2}].\] 
It follows from $m_{1}(X_{2})=m_{1}(Y_{2})$ that $\mu(A)=\mu(X)$.
\end{proof}

\begin{rem}
Note that Corollary \ref{cor-fund} also follows from the computation of $\ell^{2}$-Betti numbers of the mapping class group due to Gromov \cite{gromov-kahler} and McMullen \cite{mcmullen} (see Appendix D in \cite{kida}) and Gaboriau's work \cite{gab-l2} on the connection between $\ell^{2}$-Betti numbers and orbit equivalence. 
\end{rem}


\end{document}